\documentclass{lms}
\usepackage{amsthm, amsmath}
\usepackage{amssymb} 
\input xy
\xyoption{all}

\newtheorem{theorem}{Theorem} % 1st argument is your name for it
\newtheorem{lemma}{Lemma}    % 2nd argument is what is printed
\newtheorem*{corollary}{Corollary}

\def\Q{\mathbb{ Q}}
\def\R{\mathbb{ R}}
\def\C{\mathbb{ C}}
\def\P{\mathbb{ P}}

\def\qfl#1{\buildrel {#1}\over {\longrightarrow}}

\def\iso{\vbox{\hbox to .8cm{\hfill{$\scriptstyle\sim$}\hfill}
\nointerlineskip\hbox to .8cm{{\hfill$\longrightarrow $\hfill}} }}

\def\rond{\kern 1pt{\scriptstyle\circ}\kern 1pt}
\def\pr{\noindent\emph{Proof} : }

\newcommand\Td{\operatorname{Todd}}
\newcommand\ch{\operatorname{ch}}
\newcommand\Tr{\operatorname{Tr}}

\title[Involutions of holomorphic symplectic manifolds]% 
 {Antisymplectic involutions of holomorphic symplectic manifolds} 
\author{Arnaud Beauville}
\classno{14J50}

\begin{document}
\maketitle

\begin{abstract}
Let $X$ be a holomorphic symplectic manifold, of dimension divisible by 4, and $\sigma $ an antisymplectic involution of $X$. The fixed locus $F$ of $\sigma $ is a Lagrangian submanifold of $X $; we show that its $\hat A$-genus is $1$. As an application, we determine all possibilities for the Chern numbers of $F$ when $X$ is a deformation of the Hilbert square of a K3 surface.

\end{abstract}

\section*{Introduction}
\par  Let $X$ be an irreducible holomorphic symplectic manifold admitting an antisymplectic involution $\sigma $ (that is, $\sigma $  changes the sign of the symplectic form). The fixed locus $F$ of $\sigma $ is a Lagrangian submanifold of $X$. The main observation of this note is that \emph{when $\dim(X)$ is divisible by $4$, the $\hat A$-genus of $F$ is equal to} $1$. Our proof, given in \S 1, rests on a simple computation based on the holomorphic Lefschetz theorem.

\par  In \S 2 we apply this result when $X$ is a symplectic fourfold with $b_2=23$ (this holds when $X$ is the Hilbert square $S^{[2]}$ of a K3 surface). We show that there are exactly 11 possibilities for the pair of invariants $(K_F^2,\chi (\mathcal{O}_F))$ of the surface $F$, depending on the number of moduli of $(X,\sigma )$. In \S 3 we illustrate our results on a few examples, in particular  the \emph{double EPW-sextics} studied by O'Grady \cite{OG1}, which form the only known family of pairs $(X, \sigma )$ as above  of maximal dimension 20. 

\section{The $\hat A$-genus of the fixed manifold.}
\noindent 1.1\kern8pt Throughout  this note we consider an irreducible holomorphic symplectic 
manifold $X$ \cite{B}. This means that $X$ is compact K\"ahler, simply connected, and admits  a symplectic 2-form $\varphi \in H^0(X,\Omega ^2_X)$ which generates the $\C$-algebra $H^0(X,\Omega ^*_X)$.  We denote by $\sigma $
 an antisymplectic involution of $X$ (so that $\sigma ^*\varphi =-\varphi $).
 
\medskip
\begin{lemma}\label{lag}
The fixed locus $F$ of $\sigma $ is a smooth Lagrangian submanifold of $X$. 
\end{lemma}

\smallskip
\pr Let $x\in F$. We have a decomposition $T_x(X)=T^+\oplus T^-$ into eigenspaces  of $\sigma'(x) $. Because of the relation $\varphi _x(\sigma'(x).u,\sigma'(x).v)=-\varphi _x(u,v)$ for $u,v\in T_x(X)$, the two eigenspaces are isotropic, and therefore Lagrangian.
Since $T^+=T_x(F)$, the lemma follows.\qed
\medskip

\noindent 1.2\kern8pt Observe that the existence of the antisymplectic involution $\sigma $ forces $X$ to be \emph{projective}: indeed, let $H^2(X,\Q)^+\subset H^2(X,\R)^+$ be the $(+1)$-eigenspaces of $\sigma ^*$ in $H^2(X,\Q)\subset H^2(X,\R)$. The space $H^2(X,\R)^+$ is   contained in $H^{1,1}$, and contains a K\"ahler class; 
since $H^2(X,\Q)^+$ is dense in $H^2(X,\R)^+$, it also contains a K\"ahler class, which is ample. 

\medskip
\noindent 1.3\kern8pt   The \emph{$\hat A$-genus} $\hat A(M)$ of a compact manifold $M$ is a rational number which can be expressed as a polynomial in the Pontrjagin classes of $M$
 (\cite{H}, \S 26). When $M$ is a  complex manifold of dimension $n$, we have
\[\hat A(M)=\int_M \Td(M) \, e^{-\frac{c_1(M)}{2}} \]where $\int_M:H^*(M,\Q)\rightarrow \Q$ is the evaluation on the fundamental class of $M$
(see \cite{H}, p.~13, formula (12)). If we extend  the Euler-Poincar\'e characteristic $\chi $ as a $\Q$-linear homomorphism $K(M)\otimes \Q\rightarrow \Q$, we have $\hat A(M)=\chi (\frac{1}{2}K_M)$, where $K_M$ 
 is the canonical bundle of $M$.
 
 \medskip
\begin{theorem}
Let $X$ be an irreducible  symplectic 
manifold with $4\mid \dim(X)$, $\sigma $ an anti\-symplectic involution of $X$, $F$ its fixed manifold. Then $\hat A(F)=1$.
\end{theorem}

\smallskip
\pr Since $F$ is Lagrangian (lemma 1), the symplectic form of $X$ induces an isomorphism  $T_F\iso N^*_{F/X}$. We apply the holomorphic Lefschetz formula (\cite{AS}, 4.6):
\[\sum_i(-1)^i\Tr \sigma ^*_{\ |H^i(X,\mathcal{O}_X)}=\int_F \Td(F) (\ch \wedge N^*_{F/X})^{-1}= \int_F \Td(F) (\ch \wedge T_F)^{-1}\ .\] Because $X$ is irreducible symplectic, $\sigma ^*$ acts as $(-1)^i$ on $H^{2i}(X,\mathcal{O}_X)$; since $\dim(X)$ is divisible by 4 this implies that the above expression is equal to $1$.
\par  As usual we write the Chern polynomial $c_t(T_F)=\prod_i (1+t \gamma _i)$, where the $\gamma _i$ live in some overring of $H^*(F)$. We have 
\[\Td(F)=\prod_i\frac{\gamma _i}{1-e^{-\gamma _i}}\quad\hbox{and}\quad\ch (\wedge T_F)=\sum_{i_1<\ldots <i_k}e^{\gamma _{i_1}+\ldots +\gamma _{i_k}}=\prod_i(1+e^{\gamma _i})\ , \]hence
\[ \Td(F)(\ch \wedge T_F)^{-1}=2^{-n}e^{-c_1}\prod_i \frac{2\gamma _i}{1-e^{-2\gamma _i}}\quad,\quad\hbox{with }\ n=\dim(X)\ \hbox{ and } \ c_1=c_1(T_F)\, .
 \]
\par  Writing $\Td(F)=\sum_k \Td(F)_k$, with $\Td(F)_k\in H^{2k}(F,\Q)$, we find
\[\int_F  \Td(F)(\ch \wedge T_F)^{-1}=2^{-n}\sum_k \int_F\frac{(-c_1)^k}{k!}2^{n-k}\Td(F)_{n-k}=\int_F \Td(F) e^{-\frac{c_1}{2}}\ ,\]hence $\hat A(F)=1$.\qed

\bigskip
\par  Note that the argument applies also when $\dim(X)\equiv 2$ (mod. 4) but gives the trivial equality $\hat A(F)=0$.

\bigskip
\section{Symplectic fourfolds}
\noindent 2.1\kern8pt   When $\dim(X)=4$ the fixed locus $F$ is a surface (not necessarily connected). In that case  $\hat A(F)$ is equal to $ -\frac{1}{8}\,\mathrm{sign}(F)$, where $\mathrm{sign}(F)$ is the signature of the intersection form on $H^2(F,\R)$ (see \cite{H}, 1.5, 1.6 and 8.2.2); we have\[ \mathrm{sign}(F)=\frac{1}{3}(K_F^2-2e(F))= K_F^2-8\chi (\mathcal{O}_F)\ ,\]
where $e(F)$ is the topological Euler characteristic of $F$, and we put $K_F^2=\sum_iK^2_{F_i}$ if $F_1,\ldots ,F_p$ are the connected components of $F$.
\par  Therefore 
 Theorem 1
gives \[\mathrm{sign}(F)=K_F^2-8\chi (\mathcal{O}_F)=-8\quad \hbox{and}\quad
 K_F^2-2e(F)=-24\ .\]
\par  We will be able to say more when the action of $\sigma $ on $H^2(X)$ controls the action on $H^4(X)$, that is, when the canonical map $\mathrm{Sym}^2H^2(X)\rightarrow H^4(X)$ is an isomorphism. By \cite{G} this happens if and only if $b_2(X)=23$. This is the case for one of the two families of symplectic fourfolds known so far, namely the family of Hilbert schemes $S^{[2]}$ of  a K3 surface $S$ (and their deformations).

\medskip
\begin{theorem}
Let $X$ be a  symplectic 
fourfold with $b_2(X)=23$, $\sigma $ an antisymplectic invo\-lu\-tion of $X$, $F$ its fixed surface. Let $t$ denote the trace of $\sigma ^*$ acting on $H^{1,1}(X)$.
\par  $a)$ We have 
\[K_F^2=t^2-1\qquad \chi (\mathcal{O}_F)=\frac{1}{8}(t^2+7)\qquad e(F)=\frac{1}{2}(t^2+23) \]
\par  $b)$ The local deformation space of $(X,\sigma )$ is smooth of dimension $\frac{1}{2}(21-t)$.
\par  $c)$ The integer $t$ can take any odd value with $-19\leq t\leq 21 $. 
\end{theorem}\label{prop}

\smallskip
\pr The classical Lefschetz formula  reads
\[e(F)=\sum_i(-1)^i\Tr\, \sigma ^*_{\,|H^i(X)}\]where we put $H^*(X):=H^*(X,\Q)$.   In the case $b_2=23$, the odd degree cohomology vanishes, and the natural map $\mathrm{Sym}^2H^2(X)\rightarrow H^4(X)$ is an isomorphism \cite{G}. Let $a$ and $b$ be the dimensions of the $(+1)$- and $(-1)$-eigenspaces of $\sigma ^*$ on $H^2(X)$. We have $a+b=23$ and $a-b=t-2$. Then
\[\Tr \sigma ^*_{\,|H^4(X)}=\frac{1}{2}a(a+1)+\frac{1}{2}b(b+1)-ab=\frac{1}{2}(t-2)^2+\frac{23}{2}\ , \] 
\[e(F)=2+2\Tr \sigma ^*_{\,|H^2(X)}+\Tr \sigma ^*_{\,|H^4(X)}=2+ 2(t-2)+ \frac{1}{2}(t-2)^2+\frac{23}{2}=\frac{1}{2}(t^2+23)\ ;\leqno{\textrm{and}} \]
using  (2.1)  we deduce the other formulas  of $a)$. 
\smallskip
\par  We have $H^2(X,T_X)\cong H^2(X,\Omega ^1_X)=0$, hence the versal deformation space $\mathrm{Def}_X$ of $X$ is smooth and can be locally identified with $H^1(X,T_X)$;  the  involution $\sigma $ gives rise to an involution of  $\mathrm{Def}_X$, which under the above identification 
corresponds to $\sigma ^*$ acting on $H^1(X,T_X)$. Thus the deformation space of $(X,\sigma )$ is identified with the $(+1)$-eigenspace of $\sigma ^*$. Since $\sigma ^*\varphi =-\varphi $, this eigenspace is mapped by the isomorphism
\[ H^1(X,T_X)\qfl{i(\varphi )}H^1(X,\Omega ^1_X)\]
to the $(-1)$-eigenspace of $\sigma ^*$ in $H^1(X,\Omega ^1_X)$. With the previous notation the dimension of this eigenspace is $b-2=\frac{1}{2}(21-t)$, which proves $b)$.
\smallskip
\par Let us prove $c)$. Since $\sigma $ preserves some K\"ahler class we have  
 $a=\frac{1}{2}(t+21)\geq 1$, hence $t\geq -19$; since  $\sigma ^*\varphi =-\varphi $ we have $b=\frac{1}{2}(25-t)\geq 2$, hence $t\leq 21$. We will construct in 3.2, 3.3 and 3.4
 below examples with all possible values of $t$.\qed 

\medskip
 \begin{corollary}
The pair $(K_F^2,\chi (\mathcal{O}_F))$ can take any of the  values
$\ (0,1),\ (8,2),\ (24,4),\ (48,7),\  \allowbreak   (80,11), \   (120,16), \   (168,22), \   (224,29), \   (288, 37), \ (360, 46), \  (440, 56)$.
\end{corollary}
 
 \bigskip
 \section{Examples}
 
\noindent 3.1\kern8pt     Let $S$ be a K3 surface, $\sigma  $ an antisymplectic involution of $S$; it extends to 
 an antisymplectic involution $\sigma ^{[2]}$ of the Hilbert scheme $X=S^{[2]}$, which preserves the exceptional divisor $E$ (the locus of non-reduced subschemes).  
 We have $H^{1,1}(X)= H^{1,1}(S)\oplus \C[E]$, hence $t=\Tr  \sigma ^*_{|H^{1,1}(S)}+1$.
  The fixed locus of $\sigma $ is a curve $\Gamma $ on $S$ (not necessarily connected); the Lefschetz formula for $\sigma $ gives
 $t= e(\Gamma ) +1$. The list of all possibilities for $\Gamma $ can be found in \cite{N}.
  \par  The fixed surface $F$ of $\sigma ^{[2]}$ is the union of the symmetric square $\Gamma ^{(2)}$ and the quotient surface $S/\sigma $.
 
 \medskip
\noindent 3.2\kern8pt     Let $C$ be an irreducible plane curve of degree 6, with $s$ ordinary double points ($0\leq s\leq 10$) and no other singularities. Let $\pi :S'\rightarrow \P^2$ be
the   double covering of $\P^2$ branched along  $C$,  $S$  the minimal resolution of $S'$, and $\sigma $   the  involution of $S$ which exchanges the sheets of $\pi $. The fixed locus $\Gamma $ of $\sigma $ is the normalization of $C$; thus
 $e(\Gamma)=-18 +2s$, and $t=-17+2s$.
 
 \medskip
\noindent 3.3\kern8pt     For each integer $r$ with $1\leq r\leq 10$, there exists a K3 surface $S$ and an involution of $S$ whose fixed locus is the disjoint union of $r$ rational curves \cite{N}. Then $e(\Gamma )=2r$ and $t=2r+1$. Together with the previous example this gives all  integers $t$ appearing in  Theorem 2 $c)$, except  $t=-19$.

 \medskip
\noindent 3.4\kern8pt    The case $t=-19$ is particularly interesting, because when it holds the deformation space of $(X,\sigma )$ has maximal dimension 20 (Theorem 2, $b)$). The space  $H^2(X,\Q)^+$ is one-dimensional, generated by an ample class (1.2); the deformation space of $(X,\sigma )$ coincides locally with the  deformation space of $X$ as a polarized variety. We know only one example of this situation: O'Grady has constructed a 20-dimensional family of projective symplectic fourfolds which are double coverings of certain sextic hypersurfaces in $\P^5$, called EPW-sextics \cite{OG1}. The corresponding  involution  is antisymplectic and must satisfy $t=-19$ by Theorem 2 $b)$. The fixed surface $F$ is connected, and from Theorem 2 $a)$ we recover the invariants  $K_F^2=360$, $\chi (\mathcal{O}_F)=46$ already obtained in \cite{OG2}.

 \medskip
\noindent 3.5\kern8pt    As explained in \cite{F}\kern2pt\footnote{ I am indebted to K. O'Grady for pointing out this paper to me, thus correcting an inaccurate remark in the first version of this note.}, the above pairs $(X,\sigma )$ specialize to $(S^{[2]},\tau  )$, where $S$ is a smooth quartic surface in $\P^3$ which contains no line, and $\tau  $ associates  to a length 2 subscheme $z\in S^{[2]}$ the residual subscheme in the intersection of $S$ and the line spanned by $z$. The fixed locus becomes the surface $B$ of bitangents to $S$; this explains why $B$ has the same invariants $K_B^2=360$, $\chi (\mathcal{O}_B)=46$, as already observed by Welters \cite{W}.

 \medskip
\noindent 3.6\kern8pt    There are many other  examples, which give rise to interesting exercises. Here is one: we start with the involution $\iota $ of $\P^5$ given by $\iota  (X_0,\ldots ,X_5)=(-X_0,X_1,\ldots ,X_5)$. Let $V\subset \P^5$ be a smooth cubic threefold invariant under $\iota $: its equation must be of the form $X^2_0L(X_1,\ldots ,X_5)+G(X_1,\ldots ,X_5)=0$, where $L$ is linear and $G$ cubic. The Fano variety $X$ of lines contained in $V$ is a symplectic fourfold \cite{BD}, and $\iota $ defines an involution $\sigma $ of $X$. 
 \par  The fixed points of $\iota $ in $\P^5$ are $p=(1,0,\ldots ,0)$ and the hyperplane $H$ given by $X_0=0$. 
 A line $\ell\in X$ is preserved by $\iota  $ if and only if it contains at least two fixed points; this means that either $\ell$ contains $p$, or it is contained in $H$. The lines passing through $p$ are parametrized by the cubic surface $S\subset H$ given by $L=G=0$; the lines contained in $H$ form the Fano surface $T$ of the cubic threefold $G=0$ in $H$. Thus the fixed surface $F$ of $\sigma $ is the disjoint union of $S$ and $T$. 

\par  Using the canonical isomorphism $H^{1,1}(X)\iso H^{2,2}(V) $ \cite{BD} and Griffiths' description of the cohomology of the hypersurface $V$, one finds easily $t=-7$.
Then Theorem 2 $a)$ gives $K_F^2=48$ and $\chi (\mathcal{O}_F)=7$. Since $K_S^2=3$ and $\chi (\mathcal{O}_S)=1$, we recover the values $K_{T}^2=45$ and $\chi (\mathcal{O}_T)=6$ \cite{CG}.

\affiliationone{% in this example, two authors share an institution
   Arnaud Beauville\\
   Laboratoire J.-A. Dieudonn\'e\\
UMR 6621 du CNRS\\
Universit\'e de Nice\\
Parc Valrose\\
F-06108 Nice cedex 2, France}
  \email{arnaud.beauville@unice.fr}
% Important: Do not put any empty line here.

\bigskip






\begin{thebibliography}{99}% Replace 9 by 99 if 10 or more references
%
% Please note the use of "\and" between author names below
%

\bibitem{AS} {\bibname M. Atiyah \and I. Singer}, `The index of elliptic operators. III',  {\em Ann. of Math.} (2)  \textbf{87}  (1968), 546--604.

\bibitem{B} {\bibname A. Beauville}, `Vari\'et\'es k\"ahl\'eriennes dont la premi\`ere classe de Chern est nulle', {\em J. of Diff. Geometry} \textbf{18} (1983), 755--782. 

\bibitem{BD} {\bibname A. Beauville \and R. Donagi}, `La vari\'et\'e des droites d'une hypersurface cubique de dimension $4$', {\em C. R. Acad. Sci. Paris S\'er. I Math.}  \textbf{301}  (1985),  no. 14, 703--706. 
%\bibitem[Br]{Br} K.S. Brown: {\sl Complete Euler characteristics and fixed-point theory}.Journal of Pure and Applied Algebra \textbf{24} (1982) 103--121.

\bibitem{CG} {\bibname H. Clemens \and P. Griffiths}, `The intermediate Jacobian of the cubic threefold',
{\em Ann. of Math.} (2) \textbf{95} (1972), 281--356. 


\bibitem{F} {\bibname A. Ferretti}, `The Chow ring of double EPW sextics', {\em Preprint}	arXiv:0907.5381.
\bibitem{G} {\bibname D. Guan}, 
`On the Betti numbers of irreducible compact hyperk\"ahler manifolds of complex dimension four', {\em  Math. Res. Lett.}  \textbf{8}  (2001),  no. 5-6, 663--669. 

\bibitem{H} {\bibname F. Hirzebruch}, `Topological methods in algebraic geometry', {\em  Classics in Mathematics}. Springer-Verlag, Berlin, 1995.

\bibitem{N} {\bibname V. Nikulin}, `On the quotient groups of the automorphism groups of hyperbolic forms by the subgroups generated by $2$-reflections. Algebro-geometric applications', {\em J. Soviet Math.} \textbf{22} (1983), 1401--1476.

\bibitem{OG1} {\bibname K. O'Grady}, `Irreducible symplectic 4-folds and Eisenbud-Popescu-Walter sextics', {\em  Duke Math. J.}  \textbf{134}  (2006),  no. 1, 99--137.

\bibitem{OG2} {\bibname K. O'Grady}, `Irreducible symplectic 4-folds numerically equivalent to $(K3)^{[2]}$', {\em  Commun. Contemp. Math.}  \textbf{10}  (2008),  no. 4, 553--608.

\bibitem{W}  {\bibname G. Welters}, `Abel-Jacobi isogenies for certain types of Fano threefolds', {\em Mathematical Centre Tracts} \textbf{141}. Mathematisch Centrum, Amsterdam, 1981.\end{thebibliography}
\end{document}